\documentclass{amsart}
\usepackage{amsfonts}
\usepackage{amsthm}
\usepackage{amssymb}
\usepackage{mathrsfs}
\usepackage{epstopdf}
\usepackage{tikz}
\usepackage{tikz-cd}
\usepackage{multirow}

\def\ST{\mathrm{ST}}
\def\Proj{\mathrm{Proj}}

\def\SG{S\kern-0.02em{G}}
\def\SH{S\kern-0.02em{H}}

\def\Gal{\mathrm{Gal}}
\def\Jac{\mathrm{Jac}}
\def\Tr{\mathrm{Tr}}

\def\SL{\mathrm{SL}}

\def\PGL{\mathrm{PGL}}

\newif\iffinalrun
\iffinalrun
\else
 \fi

\iffinalrun
  \newcommand{\need}[1]{}
  \newcommand{\mar}[1]{}
\else
  \newcommand{\need}[1]{{\tiny *** #1}}
  \newcommand{\mar}[1]{\marginpar{\raggedright\tiny Fix Me:  #1 }}\fi

\usepackage{amssymb,amsmath,amsfonts,amsthm,epsfig,amscd,stmaryrd}
\usepackage{stmaryrd}
\usepackage[all,cmtip,poly]{xy}
\usepackage{url}

\newtheorem{thm2}{Theorem}
\newtheorem{cor2}[thm2]{Corollary}

\theoremstyle{definition}

\theoremstyle{remark}

\def\aae{a}
\def\abe{b}
\def\Aae{A}
\def\Abe{B}

\def\MMW{\mathcal{M}_2^w}

\newcommand{\MA}{\mathcal{A}}

\newcommand{\MM}{\mathcal{M}}
\newcommand{\MZ}{\mathcal{Z}}

\newcommand{\C}{\mathbf{C}}
\newcommand{\F}{\mathbf{F}}

\newcommand{\Q}{\mathbf{Q}}
\newcommand{\R}{\mathbf{R}}
\newcommand{\Z}{\mathbf{Z}}

\newcommand{\Qbar}{\overline{\Q}}

\newcommand\rhobar{\overline{\rho}}

\newcommand{\bbP}{\mathbf{P}}

\newcommand{\Del}{\Delta}
\newcommand{\PDel}{\Delta_{\rm poly}} 
\newcommand{\td}{\!:\!}

\DeclareMathOperator{\GL}{GL}
\DeclareMathOperator{\GSp}{GSp}
\DeclareMathOperator{\Sp}{Sp}
\DeclareMathOperator{\OO}{O}

\numberwithin{equation}{section}

\title{Abelian Surfaces with fixed~$3$-torsion}

\author[F. Calegari]{Frank Calegari}  \email{fcale@math.uchicago.edu} \address{The University of Chicago,
5734 S University Ave,
Chicago, IL 60637, USA}

\author[S. Chidambaram]{Shiva Chidambaram}  \email{shivac@uchicago.edu} \address{The University of Chicago,
5734 S University Ave,
Chicago, IL 60637, USA}

\author[D. P. Roberts]{David P. Roberts}
 \email{roberts@morris.umn.edu} \address{University of Minnesota, Morris,
600 E 4th St.,
Morris, MN 56267, USA}

\thanks{FC and SC were supported in part by NSF grant
  DMS-1701703; DPR  was supported in part by DMS-1601350. }

\newcommand{\cmmt}[1]{}

\allowdisplaybreaks
  
\begin{document}

\begin{abstract} Given a genus two curve 
$X: y^2 = x^5 + a x^3 + b x^2 + c x + d$,
we give an explicit parametrization of
all other such 
curves $Y$
with
a specified symplectic isomorphism on three-torsion of Jacobians~$\mbox{Jac}(X)[3] \cong \mbox{Jac}(Y)[3]$. 
It is known that under certain conditions
modularity of $X$ implies modularity 
of infinitely many of the $Y$, and we explain how our
 formulas render this transfer of modularity
explicit.  Our method centers on 
the invariant theory of the complex reflection
group $C_3 \times \Sp_4(\F_3)$.
We discuss other examples where complex reflection groups
are related to moduli spaces of curves, and in particular
motivate our main computation with an exposition of the simpler
case of the group~$\Sp_2(\F_3) = \SL_2(\F_3)$ and~$3$-torsion
on elliptic curves.
\end{abstract}

\maketitle

\section{Introduction}

 \subsection{Overview}
 Consider a genus two curve $X$ over $\Q$ given by an affine equation
 \begin{equation}
 \label{quintic1}
y^2 = x^5 + a x^3 + b x^2 + c x + d.
\end{equation}
The representation $\rhobar: \Gal(\overline{\Q}/\Q) \rightarrow \GSp_4(\F_3)$ on the 
three-torsion $\mbox{Jac}(X)[3]$ of its Jacobian is given by an explicit degree 80 polynomial with coefficients
in $\Q[a,b,c,d]$. 
The main theorem of this paper parametrizes all pairs $(Y,i)$ consisting 
of a curve  
\begin{equation}
\label{quintic2}
y^2 = x^5 + A x^3 + B x^2 + C x + D
\end{equation}
and a $\Gal(\overline{\Q}/\Q)$-equivariant symplectic isomorphism, 
$i : \mbox{Jac}(X)[3] \rightarrow \mbox{Jac}(Y)[3]$.  
The curves in~(\ref{quintic2}) all have a rational Weierstrass point at~$\infty$.
The reader may wonder why we did not instead try to parametrize pairs~$(Y,i)$ for~\emph{all} genus two 
curves~$Y$. The answer  is that the corresponding
moduli space, while rational over~$\C$,  will not typically be rational over~$\Q$ 
(see the discussion towards the end of \S\ref{sec:moduli}).

Analogous problems for genus one curves and 
their mod $p$ representations for $p \leq 5$ were 
solved by Rubin and Silverberg \cite{RS}.
In Section~\ref{genus1}, we explain how the mod $3$ formulas of \cite{LR} can be reconstructed by using that $\Sp_2(\F_3)$ has a two-dimensional
complex reflection representation, 
summarizing the result in Theorem~\ref{thm1}.  

Section~\ref{genus2} contains our main result, Theorem~\ref{thm2}. 
 It follows Section~\ref{genus1} closely, 
using now that $\Sp_4(\F_3)$ is the main
factor in the complex reflection group 
$C_3 \times \Sp_4(\F_3)$.     
We write the new curves as $Y = X(s,t,u,v)$ 
with $X(1,0,0,0) = X$.
The new coefficients $A$, $B$, 
$C$ and $D$ are polynomials
in $a$, $b$, $c$, $d$, $s$, $t$, $u$, and $v$.  
While the genus one and two cases are remarkably similar theoretically,
 the computations in the genus two case are  
orders of magnitude more complicated.   For example,
$A$, $B$, $C$, and $D$ have $14604$, $112763$, $515354$, and~$1727097$ terms
respectively, while the corresponding two coefficients in 
the genus one case have only $6$ and $9$ terms.  
We give all these coefficients in an accompanying 
{\em Mathematica} file available at \url{http://facultypages.morris.umn.edu/~roberts/Fixed3TorsionMma}, and this file 
also contains other material the reader may find helpful.  

Section~\ref{complements} provides four independent complements.  
\S\ref{matricial} sketches an alternative method for computing the 
above $(A,B,C,D)$. \S\ref{richelot} presents a family of  
examples involving Richelot isogenies.  \S\ref{modularity}
gives an application to modularity which was one of the motivations for this paper.
\S\ref{analogs} illustrates that much of what we do works
 for arbitrary complex reflection groups; in particular, it
sketches direct analogs of our main result in 
the computationally yet more difficult 
settings of $2$-torsion in the Jacobians of 
certain curves of genus $3$ and
$4$.

\subsection{Moduli spaces} \label{sec:moduli}  Theorems~\ref{thm1} and \ref{thm2} and the analogs sketched 
in \S\ref{analogs} are all formulated in terms of certain 
{\em a priori} complicated moduli spaces being actually open 
subvarieties of projective space.    To underscore this perspective,
we consider a whole hierarchy of standard moduli spaces
as follows.   

Let $A$ be an abelian variety over $\Q$ of dimension~$g$ with a principal polarization~$\lambda$. 
If~$V_A = A[p]$ is the set of $p$-torsion points with coefficients in 
$\overline{\Q}$, then~$V_A$ is a~$2g$-dimensional vector space 
over $\F_p$ with a symplectic form~$\wedge_A^2$ induced by the 
Weil pairing~$A[p] \times A[p] \rightarrow \mu_p$.
This structure is preserved by~$\Gal(\overline{\Q}/\Q)$, and so gives rise to a  
Galois representation: 
$$\rhobar_A:  \Gal(\overline{\Q}/\Q) \rightarrow \GSp_{2g}(\F_p);$$
here the
 similitude
 character~$\Gal(\overline{\Q}/\Q) \rightarrow \F^{\times}_p$ 
is the  mod-$p$ cyclotomic character. 

Conversely, if $\rhobar$ is any such representation on a symplectic space $(V,\wedge^2)$,
coming from an abelian variety or not,
there exists a 
moduli space~$\MA_g(\rhobar)$ 
over~$\Q$  parametrizing  triples~$(A,\lambda,\iota)$ consisting of a principally polarized abelian variety $A$ together with
an isomorphism~$\iota: (V,\wedge^2) \simeq (V_A,\wedge^2_A)$ of symplectic representations.

Via $(A,\lambda,\iota) \mapsto (A,\lambda)$, one has a covering map
$\MA_g(\rhobar) \rightarrow \MA_g$ to the moduli space of principally 
polarized $g$-dimensional abelian varieties.  For the split Galois representation~$\rhobar_0$,
 corresponding to the torsion structure
 $(\Z/p\Z)^g \oplus (\mu_p)^g$ with its natural symplectic form, 
 the cover $\MA_g(\rhobar_0)$ is the standard ``full level $p$'' cover 
 $\MA_g(p)$ of $\MA_g$.   In general, $\MA_g(\rhobar)$ is a twisted
 version of $\MA_g(p)$, meaning that the two varieties become 
 isomorphic after base change from $\Q$ to $\overline{\Q}$.

 The varieties $\MA_g(\rhobar)$ become rapidly more complicated as either 
 $g$ or $p$ increases.  In particular, they are geometrically rational exactly for 
 the cases $(g,p) = (1,2)$, $(1,3)$, $(1,5)$ $(2,2)$, $(2,3)$, and $(3,2)$ \cite[Thm II.2.1]{Hulek}. 
 In the three cases when~$g = 1$, the curves $\MA_1(\rhobar)$ are
 always rational.    In  
 the main case of interest~$(2,3)$ for this paper,  the three-dimensional variety $\MA_2(3) = \MA_2(\rhobar_0)$ is rational~\cite{Bruin}.
However, for many~$\rhobar$, 
including all surjective representations,
it is proven in~\cite{CC} that the variety $\MA_2(\rhobar)$
 is never rational.
 It \emph{is} true, however, that
 there exists a degree~$6$ cover~$\MA^{w}_2(\rhobar)$ which is rational  (\cite[Lemma~10.2.4]{BCGP}).
 Thus while Theorem~\ref{thm1}  corresponds to a parametrization of~$\MA_1(\rhobar)$
for~$p = 3$, 
  Theorem~\ref{thm2} corresponds to a parametrization 
 of~$\MA^{w}_2(\rhobar)$. 
 More precisely, the Torelli map~$\MM_2 \rightarrow \MA_2$ is an open immersion,
 and the pullback of~$\MA^{w}_2(\rhobar)$ is 
 the moduli space~$\MMW(\rhobar)$ of genus two curves
of the form~\eqref{quintic1}  whose Jacobians give rise to~$\rhobar$,
 and it is~$\MM^{w}_2(\rhobar)$ which we explicitly parametrize.
The retreat to this cover
 is optimal in the sense that six is  
generically the minimal degree of any dominant
rational map from~$\mathbf{P}^3_{\Q}$ to~$\MA_2(\rhobar)$ \cite{CC}.  
We mention in passing that our arguments give an alternative proof of~\cite[Lemma~10.2.4]{BCGP}.

  There is a natural generalization of the varieties $\MA_g(\rhobar)$.  
      Namely, for any $m \in \F_p^\times$, 
 one can require instead an isomorphism $i : (V,\wedge^2) \simeq (V_A,m \wedge^2_A)$.  
 For $m/m'$ a square, the corresponding
 varieties are canonically isomorphic, so that one 
 gets something new only in the case of $p$ odd.
 We denote this new moduli space involving ``antisymplectic'' isomorphisms
  by $\MA_g^*(\rhobar)$.   
 Our policy throughout this paper is to focus on 
 $\MA_g(\rhobar)$ and be much briefer about
 parallel results for $\MA_g^*(\rhobar)$.

\section{Elliptic curves with fixed $3$-torsion}
\label{genus1}
   In this section, as a warm up to Section~\ref{genus2}, 
   we rederive 
   formulas describing elliptic curves with fixed
3-torsion from the invariant theory of the group $\Sp_2(\F_3)$.   Many of the
steps in the derivation transfer with no theoretical change to our main case
of abelian surfaces.   We present these steps in greater detail here, 
because space allows us to give explicit formulas right in the text.  Throughout this 
section and the next, we present the derivations in elementary language
which stays very close to the computations involved.   Only
towards the end of the sections do we recast the results in the
moduli language of the introduction.

\subsection{Elliptic curves and their $3$-torsion}
     Let $a$ and $b$ be rational numbers such that the polynomial discriminant
 $\PDel = -4 a^3 - 27 b^2$ of $x^3+ax+b$
 is nonzero and consider the elliptic curve $X$ over $\Q$ with affine equation 
 \begin{equation}
 \label{elliptic}
 y^2=x^3+ax+b.
 \end{equation} 
 We emphasize the discriminant $\Del(a,b) = \Del = 2^4 \PDel$ 
 in the sequel,  because it makes  \S\ref{finding1} cleaner.
 
 By a classical division polynomial
 formula, the eight primitive $3$-torsion points $(x,y) \in \C^2$
 are exactly the points satisfying both \eqref{elliptic} and 
 \begin{equation}
 \label{elliptic3}
 3 x^4 + 6 a x^2 + 12 b x -a^2 = 0.
 \end{equation}
 Equations~\eqref{elliptic} and \eqref{elliptic3} together define an 
 octic algebra over $\Q$.  Rather than work with the two generators 
 $x$ and $y$ and the two relations \eqref{elliptic} and \eqref{elliptic3}, 
 we shift over to a more standard number-theoretic
 context by eliminating $x$ via a resultant.  Then the algebra
 in question is the quotient $K:=K_{a,b}$ of $\Q[y]$ coming from the equation
 \begin{equation}
 \label{division1}
  2^8 3^3 y^8 + 2^{11} 3^3 b y^6 - 2^5 3^2 b \Delta y^4 - \Delta^2 = 0.
 \end{equation}
 The left side of \eqref{division1} is the classical division polynomial associated
 to the elliptic curve \eqref{elliptic}.  
 
 It will be convenient to work with a variant of \eqref{division1} 
 where $a$ plays a more prominent role.   It turns out that 
 $3x \in K$ has square roots in $K$.  We pick one of them 
 and call it $z$ by the formula
  \[
 z=-24 y \Delta^{-1} \left(10 a^2 x+9 a b+6 a x^3-9 b x^2\right).
 \]
 Then, assuming $a \neq 0$ to avoid inseparability issues,
  $K$ is given as the quotient of $\Q[z]$ coming from the equation
 \begin{equation}
 \label{division3}
 F(a,b,z) := z^8 + 18 a z^4 + 108 b z^2 - 27 a^2 = 0.
 \end{equation}
 We will use this division polynomial rather than \eqref{division1}
 in the sequel.

\subsection{$\Sp_2(\F_3)$ and related groups}
For generic $(a,b)$, the Galois group of the polynomial $F(a,b,z)$ is 
$\GSp_2(\F_3) = \GL_2(\F_3)$.   The discriminant of $F(a,b,z)$ is 
$-2^{8} 3^{21} a^2 \Delta^4$.   Thus the splitting
field $K'_{a,b}$ of $F(a,b,z)$ contains $E = \Q(\sqrt{-3})$ 
for all $a$, $b$.  The relative 
Galois group $\mbox{Gal}(K'_{a,b}/E)$ is 
$\Sp_2(\F_3) = \SL_2(\F_3)$.   
We will generally use symplectic  
rather than linear language in the sequel, to
harmonize our notation with our main case of
genus two.   
Also we will systematically use $\omega = \exp(2 \pi i/3) = (-1+\sqrt{-3})/2$
as our preferred generator for $E$.  

To describe elliptic curves with fixed $3$-torsion, 
we use that \eqref{division3} arises as a generic polynomial in the invariant 
theory of $\Sp_2(\F_3)$.   The invariant theory is 
simple because  ${\renewcommand{\arraycolsep}{1pt} \Sp_2(\F_3) = \langle  \mbox{\scriptsize{$ \left( \begin{array}{cc} 1 & 0 \\ 1 & 1 \end{array} \right),  \left( \begin{array}{cc} 1 & 1 \\ 0 & 1\end{array} \right)$}}  \rangle}$ can be realized as a complex reflection group by sending the generators in
order to 
 \begin{align}
 \label{mat2}
 g_1 & = \left( \begin{array}{cc}  
  \bar{\omega} &    \bar{\omega}-1  \\
 0 & 1 \\
                    \end{array} \right), &
g_2 & =  \left(
                 \begin{array}{cc}
                  1 & 0  \\
                   (\omega-1)/3 & \omega \\
                 \end{array}
                 \right).
\end{align}
The matrices $g_1$ and $g_2$ are indeed complex reflections because all 
but one eigenvalue is $1$.  In our study of the image $\ST4 = G = \langle g_1,g_2 \rangle$,
the subgroup  $H = \langle g_1 \rangle$ will play an important role.   
Here our notation $\ST4$ refers to the placement of $G$ in the
Shephard--Todd classification of the thirty-seven exceptional irreducible complex reflection
groups sorted roughly by increasing size \cite[Table VII]{ShTo}.  

For both the current case of $n=2$ and the main case of $n=4$, we 
are focused principally on three irreducible characters of 
$\Sp_n(\F_3)$, the unital character $\chi_1$ and a complex conjugate
pair $\chi_{na}$ and $\chi_{nb}$.  Here 
$\chi_{na}$ corresponds to the representations \eqref{mat2} and \eqref{mat4} on
$V = E^n$.   Just as {\em invariant} is used for polynomials associated to $\chi_1$, 
we will use the terms {\em covariant} and {\em contravariant}
for polynomials similarly associated to $\chi_{na}$ and $\chi_{nb}$ 
respectively.  

The left half of Table~\ref{Sp2F3} shows how the three characters
$1$, $\chi_{2a}$, and $\chi_{2b}$ fit into the entire character
theory of $\Sp_2(\F_3)$.   For example, via $\bar{\omega}+1=-\omega$ and its
conjugate, $g_1$ and $g_2$ lie in the classes $3A$ and $3B$ respectively.  
While this information is clarifying,
it is not strictly speaking needed for our arguments.

\begin{table}[htb]
\[
{\renewcommand{\arraycolsep}{4pt}
 \begin{array}{c|rrrrrrr|ccccccccc| l }
        |C|:   &   1  &  1  & 4 & 4 & 6 & 4 & 4 &  \multicolumn{9}{c|}{\langle \chi_i,\phi_k \rangle} &     \\
 C:&1A & 2A & 3A & 3B & 4A & 6A & 6B & 0 & 1 & 2 & 3 & 4 & 5 & 6 & 7 & 8 & N_i(x) \\
 \hline
    \chi_1 & 1 & 1 & 1 & 1 & 1 & 1 & 1 & 1 & \text{} & \text{} & \text{} & 1 &
      \text{} & 1 & \text{} & 1 & 1 \\
    \chi_{1a} & 1 & 1 & \bar{\omega} & \omega  & 1 &  \bar{\omega} & \omega  &
      \text{} & \text{} & \text{} & \text{} & 1 & \text{} & \text{} & \text{} &
      1 & x^4 \\
    \chi_{1b} & 1 & 1 & \omega  &  \bar{\omega} & 1 & \omega  &  \bar{\omega}&
      \text{} & \text{} & \text{} & \text{} & \text{} & \text{} & \text{} &
      \text{} & 1 & x^8 \\
    \chi_2& 2 & -2 & -1 & -1 & 0 & 1 & 1 & \text{} & \text{} & \text{} & \text{} &
      \text{} & 1 & \text{} & 1 & \text{} & x^5+x^7 \\
    \chi_{2a} & 2 & -2 & -\omega  &  -\bar{\omega} & 0 & \omega  &  \bar{\omega}&
      \text{} & 1 & \text{} & 1 & \text{} & 1 & \text{} & 2 & \text{} & x+x^3 \\
    \chi_{2b} & 2 & -2 & - \bar{\omega}& -\omega  & 0 &  \bar{\omega} & \omega  &
      \text{} & \text{} & \text{} & 1 & \text{} & 1 & \text{} & 1 & \text{} &
      x^3+x^5 \\
    \chi_3 & 3 & 3 & 0 & 0 & -1 & 0 & 0 & \text{} & \text{} & 1 & \text{} & 1 &
      \text{} & 2 & \text{} & 2 & x^2+x^4+x^6 \\
   \end{array}
   }
 \]
 \caption{\label{Sp2F3} Character table of $\Sp_2(\F_3)$ and invariant-theoretic information}
 \end{table}
 
The right half of Table~\ref{Sp2F3} gives numerical information that will guide our
calculation with explicit polynomials in the next subsections.  
 The characters
are orthonormal with respect to the Hermitian inner product 
$\langle f,g \rangle = |G|^{-1} \sum_{C} |C| f(C) \overline{g(C)}$.
 Let $\phi_k = \sum_i \langle \chi_i, \phi_k \rangle \chi_i$ be the character of the $k^{\rm th}$ symmetric
 power $\mbox{Sym}^k V$.   The multiplicities $\langle \chi_i, \phi_k \rangle$
 for $k \leq 8$ are given in the right half of 
 Table~\ref{Sp2F3}.   These numbers are given for arbitrary $k$ by 
 $
 \sum_{k=0}^\infty  \langle \chi_i, \phi_k \rangle x^k = N_i(x)/((1-x^4)(1-x^6)).
 $
 The character of the permutation representation of $G$ on the coset space $G/H$ is
 $\phi_{G/H} = \chi_1 + \chi_3 + \chi_{2a} + \chi_{2b}$.  
 If $W$ has character $\chi_i$ then the dimension of the subspace $W^H$ of $H$-invariants is 
  $\langle \chi_i, \phi_{G/H} \rangle$.
 So $\dim(W^H)=1$ if $i \in \{1,2a,2b,3\}$ and
 $\dim(W^H)=0$ if $i \in \{1a,1b,2\}$.    
 
\subsection{Rings of invariants}
\label{section:invariants1}
The group $G$ acts on the polynomial ring $E[u,z]$ by the formulas induced from the matrices in \eqref{mat2},
\begin{align*}
g_1 u & = \bar{\omega}u + (\bar{\omega}-1)z, & g_2 u & = u,  \\ 
g_1 z & = z, & g_2 z & = (\omega-1)u/3 + \omega z . 
\end{align*}
Despite the appearance of the irrationality $\omega$ in these formulas, 
there is an important rationality present.  Namely we have arranged
in \eqref{mat2} that  $g_1^2 = \overline{g}_1$ and $g_2^2 = \overline{g}_2$.
Accordingly $G$ is stable under complex conjugation,
a stability not present in either the original Shephard--Todd paper \cite[\S4]{ShTo}
or in {\em Magma}'s implementation \verb@ShephardTodd(4)@.

We can use stability under 
complex conjugation to interpret $G$ and $H$ as the 
$E$-points of  group schemes $\underline{G}$ and $\underline{H}$ over $\Q$.   
Then actually $\underline{G}$ acts on 
$\Q[u,z]$.  All seven irreducible representations of $\underline{G}$ are defined over
$\Q$, just like all three representations of the familiar group scheme $\underline{H} \cong \mu_3$,
 are defined over
$\Q$.   In practice, we continue thinking almost exclusively in terms of ordinary groups;
these group schemes just provide a conceptually clean way of saying that in 
our various choices below we can and do always take all coefficients rational.  

Define
\begin{align}
\label{wab}
w & = \frac{u^3}{3} + u^2 z + u z^2, & 
a & =  \frac{w z}{9} , & 
b & = \frac{w^2-6 w z^3-3 z^6}{324}
\end{align}
in $\Q[u,z]$.  Then the subrings of $\underline{H}$- and $\underline{G}$-invariants are
respectively
\begin{align}
\label{invariants1}
\Q[u,z]^{\underline{H}} &= \Q[w,z], &
\Q[u,z]^{\underline{G}} & = \Q[a,b].
\end{align}
Giving $u$ and $z$ weight one, the elements $w$, $a$, and $b$ clearly have weights $3$, $4$, and $6$ 
respectively. 
If one eliminates $w$ from the last two equations of \eqref{wab}, then one gets
the polynomial relation $F(a,b,z)=0$ of \eqref{division3}, explaining our choice
of overall scale factors in \eqref{wab}.   The fact that the rings on the right
in \eqref{invariants1} are polynomial rings, rather than more complicated
rings requiring relations to describe, comes exactly from the fact that
$H$ and $G$ are complex reflection groups, by the Chevalley--Shephard--Todd theorem~\cite{CST}. 

\subsection{Covariants and contravariants} 
The graded ring $\Q[w,z]$ is free of rank eight over the graded ring 
$\Q[a,b]$.   Moreover there is a homogeneous basis
$1$, $z^2$, $z^4$, $z^6$, $\alpha_1$, $\alpha_3$, $\beta_3$, $\beta_5$ 
with the following properties.   The exponent or index $d$ gives the 
weight, and the elements $\alpha_d$ and 
$\beta_d$  are in the isotypical piece of $\Q[u,z]_d$ corresponding to 
$\chi_{2a}$ and  $\chi_{2b}$ respectively.   

The {covariants} $\alpha_d$ and the {contravariants} $\beta_d$
are each well-defined up to multiplication by a nonzero rational scalar.  
Explicit formulas for particular choices can be found by simultaneously 
imposing the $G$-equivariance condition and 
the $H$-invariance condition.  We take
\begin{align}
\label{alphabeta} \alpha_1 & = z , &   
\alpha_3 & = \frac{w+z^3}{6} ,&  \beta_3  & = \frac{w-z^3}{2}, &  \beta_5 & = \frac{5 w z^2+3 z^5}{18}.
\end{align}
Ideas from classical invariant theory are useful in finding these
quantities.   For example, the polynomials in $\Q[u,z]_3$
which have the required $G$-equivariance property for contravariance are
exactly the linear combinations of the partial derivatives
$\partial_u a$ and $\partial_z a$.  The subspace fixed 
by $H$ is the line spanned by
$(\partial_u - \partial_z) a$.  Thus $\beta_3 \propto (\partial_u - \partial_z) a$ and,
in the same way, $\beta_5 \propto (\partial_u -  \partial_z) b$.

\subsection{New coefficients}  
\label{charpol1} \label{newcoeff1} The octic $\Q[a,b]$-algebra $\Q[w,z]$ acts on 
itself by multiplication and so every element $e$ in $\Q[w,z]$ has an octic characteristic
polynomial $\phi(e,u) \in \Q[a,b,u]$.   One has $\phi(z,u) = F(a,b,u)$ from \eqref{division3}. 
To obtain the characteristic polynomial for a general $e$, one can express $e$ as an 
element of $\Q(a,b,z)$ via \eqref{alphabeta} and $w = 9a/z$.   Then 
one removes $z$ by a resultant to get the desired octic relation     
on $e$.   Alternatively, we could have calculated these characteristic polynomials
by using 8-by-8 matrices; in \S\ref{newcoeffs2} we use the matrix approach.  

Carrying out this procedure for the general covariant and
contravariant gives
\begin{align*}
\phi(s \alpha_1+t \alpha_3,u) & = F(A(a,b,s,t),B(a,b,s,t),u), \\
\phi(s \beta_3+t \beta_5,u) & = F(A^*(a,b,s,t),B^*(a,b,s,t),u), 
\end{align*}
with 
\begin{eqnarray*}
	3A(a,b,s,t) & = & 3 a s^4 +18 b s^3 t -6 a^2 s^2 t^2 -6 a b s t^3 -(a^3+9
	b^2) t^4, \\
	9B(a,b,s,t) & = & 9 b s^6-12 a^2 s^5 t-45 a b s^4 t^2-90 b^2 s^3 t^3 + 15 a^2 b s^2 t^4  \\
	&& \qquad -2 a 
	  (2 a^3+9 b^2 ) s t^5 -3 b   (a^3+6 b^2 ) t^6,
\end{eqnarray*}
and $A^*$ and $B^*$ in the accompanying computer file. 
As stated in the introduction, $A$ and $B$ when fully expanded have $6$ and
$9$ terms respectively and agree exactly with expressions in \cite[\S2]{LR}.

The polynomials $A$ and $B$ and their starred versions are
 respectively of degrees four and six in $s$
and $t$.  Also in the main case assign weights $(4,6,-1,-3)$ to 
$(a,b,s,t)$ and in the starred case make these weights $(4,6,-3,-5)$ instead.  
Then all four polynomials are homogeneous of weight zero.   

\subsection{Geometric summary}
\label{geom1} The following theorem summarizes
our calculations in terms of moduli spaces.   The $\rhobar$ of the introduction
is the mod $3$ representation of the initial elliptic curve, so to be more
explicit we write $\MA_{a,b}$ rather than $\MA_1(\rhobar)$.  
\begin{thm2}
\label{thm1}
Fix an equation $y^2=x^3+ax+b$ defining an elliptic curve $X$ over $\Q$. 
Let $\MA_{a,b}$ be the moduli space 
of pairs $(Y,i)$ with $Y$ an elliptic curve and $i : X[3] \rightarrow Y[3]$ 
a symplectic isomorphism.  Then $\MA_{a,b}$ can be realized
as the complement of a discriminant locus $\MZ_{a,b}$ in the projective 
line $\Proj \; \Q[s,t]$.   The 
natural map to the $j$-line ${\MA}_1 \subset \Proj \; \Q[A,B]$ 
has degree twelve and is given by 
\begin{equation}
\label{maineq1}
(A,B) = (A(a,b,s,t),B(a,b,s,t)).
\end{equation}
The formula~$y^2 = x^3 + A(a,b,s,t) x + B(a,b,s,t)$ gives the universal elliptic curve $X(s,t)$ over $\MA_{a,b}$. 
\end{thm2}
\noindent The discriminant locus $\MZ_{a,b}$ is
given by the vanishing of the discriminant 
\begin{equation}
\label{disc1}
\Delta(A(a,b,s,t),B(a,b,s,t)) = 27 \Delta(a,b) (3 t^4 + 6 a s^2 t^2 + 12 b s t - a^2 s^4)^3.
\end{equation}
It thus consists of four geometric points.  Comparing with \eqref{elliptic3}, one 
sees that these points are permuted by $\Gal(\overline{\Q}/\Q)$ 
according to the projective mod $3$ representation into 
$\PGL_2(\F_3) \cong S_4$.   Theorem~\ref{thm1} has 
a direct analog for the covers $\MA_{a,b}^* \rightarrow \MA_1$.

\subsection{Finding $(s,t)$}  
\label{finding1} Let $X: y^2 = x^3+a x +b$ and 
$Y: y^2 = x^3 + Ax + B$ be elliptic curves over $\Q$ 
with isomorphic $3$-torsion.  Then, in contrast with 
the analogous situation for the 
genus two case described in \S\ref{finding2}, it is very easy
to find associated $(s,t) \in \Q^2$.  Namely, \eqref{maineq1} 
and its analog $(A,B) = (A^*(a,b,s,t),B^*(a,b,s,t))$ each have
twenty-four solutions in $\C^2$.    One just
extracts the rational ones, say by eliminating 
$s$ and factoring the resulting degree twenty-four polynomials 
$f(t)$ and $f^*(t)$.   If the image of $\Gal(\overline{\Q}/\Q)$
is all of $\GSp_2(\F_3) = \GL_2(\F_3)$, then one of
these polynomials factors as $1+1+6+8+8$ and
the other as $12+12$.  The two $1$'s correspond 
to the desired solutions $\pm(s,t)$.

Discriminants are useful in distinguishing the two moduli spaces as follows. 
If $Y$ has the form $X(s,t)$ then $\Delta_X/\Delta_Y$ is 
a perfect cube by \eqref{disc1}.  If it has the form $X^*(s,t)$ then $\Delta_X \Delta_Y$
is a perfect cube by the starred analog of \eqref{disc1}.  These implications determine a unique space on 
which $Y$ represents a point unless  $\Delta_X$ and $\Delta_Y$ 
are both perfect cubes.   Since $x^3-\Delta$ is a resolvent
cubic of the octic \eqref{division3}, this ambiguous
case arises if and only if the image $\Gamma$ of $\rhobar_X$ has
order dividing 16.  

As an example, let $(a,b) = (-1,0)$ so that 
$X$ has conductor $2^5$ and discriminant $2^6$.
Let $(A,B) = (-27,-162)$ so that $Y$ has conductor $2^5 3^3$ 
and discriminant $-2^9 3^9$.   The octic polynomials 
$F(a,b,z)$ and $F(A,B,z)$ define the 
same field because under {\em Pari}'s 
\verb@polredabs@ they each become 
$z^8+6 z^4 - 3$.  This polynomial
has Galois group of order $16$.  
The procedure in the first paragraph
yields solutions only in the starred
case, these being $(s,t) = \pm (-1/2,3/2)$.
A Galois-theoretic argument shows 
that one gets solutions in both cases
if and only if $\Gamma$ is at most a quadratic 
 twist of
the representation~$\rhobar_0 = \Z/3\Z \oplus \mu_3$.
 An instance 
is $X=Y$ coming
from $(a,b) = (5805, -285714)$ which is the modular curve~$X_0(14)$ of genus one 
and discriminant $-2^{18} 3^{12} 7^3$;
here $(s,t)=\pm(1,0)$ in the main case and $2^6 3^4 7^2 (s,t) = \pm(435,11)$ in 
the starred case.

\section{Abelian surfaces with fixed $3$-torsion}
\label{genustwo} \label{genus2}
   In this section, we present our main theorem on abelian surfaces with 
fixed $3$-torsion.  We are brief on parts of the derivation which closely follow steps 
described in the previous section, and concentrate on steps which have a new
feature.

\subsection{Weierstrass curves and their $3$-torsion}  
\label{3tors2}
By a {\em Weierstrass curve} in 
this paper we will mean a genus two curve together with a distinguished Weierstrass point.  
Placing this marked point at infinity and shifting the variable $x$,
one can always present a Weierstrass curve via the affine equation  
\eqref{quintic1}, which we call a~{\em Weierstrass equation}. 
Replacing $(a,b,c,d)$ by $(u^4 a, u^6 b, u^8 c, u^{10} d)$ yields an isomorphic
Weierstrass curve via the compensating change $(x,y) \mapsto (u^2 x,u^5 y)$.
The standard discriminant of the genus two curve \eqref{quintic1} is $\Del(a,b,c,d) = \Del = 2^8 \PDel$,
where $\PDel$ is the discriminant of the quintic polynomial on the right of
\eqref{quintic1}.    
It is best for our purposes to give the parameters
$a$, $b$, $c$, and $d$ weights $12$, $18$, $24$, and $30$.  
In this system, $\Delta$ is homogeneous of weight $120$. 
The (coarse) moduli space of Weierstrass curves $\MMW$ is then  
the complement of the hypersurface $\Delta=0$ in the weighted
projective space $\bbP^3(12,18,24,30) = \bbP^{3}(2,3,4,5)$.
As explained at the end of \S1.2, rather than describing moduli spaces mapping to $\MA_2$, we will be describing 
their base changes to $\MM_2^w$.

The group law in terms of effective divisors on the Jacobian of a general genus two curve $X: y^2 = f(x)$ 
yields a classical $\Gal(\overline{\Q}/\Q)$-equivariant bijection \cite{DD} from the non-zero $3$-torsion points to 
decompositions of the form
\[
 f(x) = (b_4 x^3 + b_3 x^2 + b_2 x + b_1)^2 - b_7 (x^2 + b_6 x + b_5)^3.
 \]
 In the quintic case of  \eqref{quintic1},  one has $b_4^2 = b_7$. The minimal polynomial of $b_4^{-2}$ is a degree $40$ polynomial $p_{40}$ such that $p_{40}(x^2)$ describes the $3$-torsion representation
 of $X$.

 In our reflection group approach, it is actually $p_{40}(z^6)$ which
 appears naturally.   It has $1673$ terms and begins as
\begin{eqnarray}
\nonumber \lefteqn{F(a,b,c,d,z) =  z^{240} +15120 a z^{228} +  2620800 b z^{222}  } &&  \\ 
\label{deg240} &&   -504  \left(70227 a^2-831820 c\right) z^{216}  
-1965600 z^{210} (2529 a b-33550
d) z^{210}  + \cdots .
\end{eqnarray}
The splitting field of $F(a,b,c,d,z)$ is the compositum of the splitting fields of 
$p_{40}(x^2)$ and $x^3-\Delta$.  
In particular, having chosen a \emph{Weierstrass equation},
the field~$E(\Delta^{1/3})$ remains constant throughout  our 
family of Weierstrass equations, even though~$E(\Delta^{1/3})$  is \emph{not} determined by the~$3$-torsion representation.
On the other hand, the change of coordinates
$(x,y) \mapsto (u^2 x, u^5 y)$ maps~$\Delta$ to~$u^{40} \Delta$, and so this auxiliary
choice places no restrictions on the \emph{Weierstrass curves} which can occur in the family.
In contrast, when~$g=1$, the field~$E(\Delta^{1/3})$   also remains constant, but
in this case it \emph{is} determined by the~$3$-torsion representation as it is the 
 fixed field of the~$2$-Sylow of the image of~$\mbox{Gal}(\overline{\Q}/E)$ in~$\Sp_2(\F_3)$.

\subsection{$\Sp_4(\F_3)$ and related groups} 
Define $g_1$, $g_2$, $g_3$, and $g_4$ to be
\begin{equation}
\label{mat4}
{\renewcommand{\arraycolsep}{3pt}
           \left(
                     \begin{array}{cccc}
                      1 & 0 & 0 & 0 \\
                      0 &1 & 0 & 0 \\
                      0 & 0 & \omega & 0 \\
                      0 & 0 & 0 & 1 \\
                     \end{array}
                     \right),  \;
                     \left(
                     \begin{array}{cccc}
                      \alpha & -\bar{\alpha} & -\bar{\alpha} & 0 \\
                      -\bar{\alpha} & \alpha & -\bar{\alpha} & 0 \\
                      -\bar{\alpha} & -\bar{\alpha} & \alpha & 0 \\
                      0 & 0 & 0 & 1 \\
                     \end{array}
                     \right), \;
  \left(
                     \begin{array}{cccc}
                      1 & 0 & 0 & 0 \\
                      0 & \omega & 0 & 0 \\
                      0 & 0 & 1 & 0 \\
                      0 & 0 & 0 & 1 \\
                     \end{array}
                     \right), \;
                     \left(
                     \begin{array}{cccc}
                      \alpha & \bar{\alpha} &0& \bar\alpha \\
                       \bar{\alpha} & \alpha &0& -\bar{\alpha} \\
                       0&0&1&0 \\
                      \bar{\alpha} & -\bar{\alpha} &0 & \alpha \\
                     \end{array}
                     \right),
}
\end{equation}
where $\alpha = \omega/\sqrt{-3}$.  Define  $H = \langle g_1,g_2,g_3 \rangle$  and
 $G = \langle g_1,g_2,g_3,g_4 \rangle$.   The matrices $g_i$ are all complex reflections of
 order $3$, and they are exactly the matrices given in \cite[10.5]{ShTo}.  
   As with~$H = C_3$ and~$G=\ST4=\Sp_2(\F_3)$ of the last section, the new groups~$H=\ST25$ 
and~$G=\ST32$ are also complex reflection groups.   The group $G$ has the structure  
$C_3 \times \Sp_4(\F_3)$ and it is the extra $C_3$ that is the reason
that $\Delta$ behaves differently in the two cases.

  Again numeric identities guide polynomial calculations
 as we discussed around Table~\ref{Sp2F3}.  For example, orders are products
 of degrees of fundamental invariants.  Analogous to the old cases $|C_3|=3$ and
 $|\Sp_2(\F_3)| = 4 \cdot 6$, the new cases are $|H|=6 \cdot 9 \cdot 12$ and 
 $|G|=12 \cdot 18 \cdot 24 \cdot 30$.   Thus again the index $|G|/|H|=240$ 
 matches the degree of the main polynomial \eqref{deg240}.  
 The character table of $G$ has size $102 \times 102$, 
 so we certainly will not present the analog of Table~\ref{Sp2F3}.  
 The most important information is that the degrees in which 
 co- and contravariants live, previously 1, 3 and 3, 5, are now
 $1$, $7$, $13$, $19$ and $11$, $17$, $23$,  $29$ for $G$.   
 
 \subsection{Rings of invariants}
 \label{invariants2}
 One has the rationality condition $g_i^2 = \overline{g}_i$ for
 all  four $i$, allowing us again to interpret $H$ and $G$ as $E$-points of
 group schemes $\underline{H}$ and $\underline{G}$ over $\Q$.  
 The matrices $g_i$ together give an action of $\underline{G}$ on 
 $\Q[z_1,z_2,z_3,z_4]$.   The variable $z=z_4$ plays a role which
is different from the other $z_i$.

Define, following \cite[4.72]{Hunt},
  \begin{eqnarray*}
  p&=&z_1^6 + z_2^6 +z_3^6 -10 \left( z_2^3 z_3^3+z_2^3 z_1^3+z_3^3 z_1^3\right), \\
  q & =  &  (z_1^3-z_2^3)(z_1^3-z_3^3)(z_2^3-z_3^3), \\
  r & = & (z_1^3 + z_2^3 +z_3^3) \left[ (z_1^3 + z_2^3 +z_3^3)^3 + 216 z_1^3 z_2^3 z_3^3 \right].
  \end{eqnarray*}
Define also $a$, $b$, $c$, and $d$ by 
taking~$2^4 3^7 5 a$, $2^6 3^9 5^2 b$, $2^8 3^{12} 5^3 c$, $2^{10} 3^{16} 5^5 d$ to be
  {\small
 \[
 \begin{array}{c}
 - p^2-5  r+1320  q z^3-132  p z^6-6 z^{12},  \\[.04in]
 p^3-400  q^2-5  p  r-680  p  q z^3+ 323
     p^2 z^6 -255  r z^6-7480  q z^9+68  p z^{12}-4 z^{18}, \\[.04in]
2  p^4-800  p  q^2-5  p^2  r+ 320  p^2   q z^3 -3000  q
     r  z^3 -722  p^3 z^6 +175200  q^2 z^6 +990  p  r  z^6 \\  
 +33040  p  q z^9 -953  p^2 z^{12} +3495  r 
    z^{12}+15720  q z^{15}+268  p z^{18}-3 z^{24}, \\[.04in]
13  p^5-6000  p^2  q^2-25  p^3  r+ 21600  p^3  q z^3 -9600000
     q^3 z^3 -45000  p  q  r  z^3 
      + 11790  p^4 z^6 \\ -4572000  p  q^2 z^6 -37575  p^2  r z^6 +28125  r^2 
    z^6 - 247200  p^2  q z^9 +945000  q  r  z^9 \\  + 37155  p^3 z^{12} +234000  q^2 z^{12} 
    -150075  p  r     z^{12}-214200  p  q z^{15}+ 30855  p^2 z^{18} \\ -143775  r  z^{18} 
     +354600  q z^{21}+2340  p z^{24}-12
    z^{30}.
\end{array}
\]
}%
 Because $H$ and $G$ are complex
  reflection groups, the rings of invariants are freely generated, 
 explicit formulas being 
  \begin{align*}  
\label{invariant2} 
 \Q[z_1,z_2,z_3,z]^{\underline{H}} & = \Q[ p, q, r,z], & 
 \Q[z_1,z_2,z_3,z]^{\underline{G}} & = \Q[a,b,c,d].  
 \end{align*}
 When one removes $p$, $q$, $r$ from the equations defining $a$, $b$, $c$, $d$, one gets exactly the 
 degree $240$ equation \eqref{deg240} for $z$.

 \subsection{Covariants and contravariants} 
 \label{covariants2}
 As mentioned before, group-theoretic calculations like those in Table~\ref{Sp2F3} say that covariants
 lie in degrees $1$, $7$, $13$, and $19$.   Formulas for $H$-invariant covariants
 in these degrees are 
 {\small
  \begin{align*}
 \alpha_1 & = z,  & 
 2^2 3^3 5 \alpha_7 & = 7  p z-3 z^7,  &
 2^4 3^6  \alpha_{13} & =  (11  r-3  p^2)z+216  q z^4 +72  p z^7,  \\
 2^3 3^{10}   \alpha_{19}  &   \omit\rlap{$\;=  ( p^3- p  r-468  q^2)z -24  p  q z^4 +z^7 (66  r-6  p^2)-288  q z^{10} -12  p  z^{13}.$} 
 \end{align*}
}%
  Here, unlike in the genus one case, there is an ambiguity beyond multiplying by a nonzero scalar.
 Namely rather than working with $\alpha_{13}$ we could work with any linear
 combination of $a \alpha_1 $ and $\alpha_{13}$ that involves $\alpha_{13}$ nontrivially.
 Similarly we could replace $\alpha_{19}$ by
  $c_1 b \alpha_1  + c_7 a \alpha_7  + c_{19} \alpha_{19}$ for any nonzero $c_{19}$.
The choices involved in picking particular contravariants $\beta_k$ mirror the
choices involved in picking $\alpha_{k-10}$.  Our choice of 
 $(\beta_{11},\beta_{17},\beta_{23},\beta_{29})$ is given in the accompanying computer
file.  

\subsection{New coefficients} \label{newcoeff2}
\label{newcoeffs2} Our key computation  
takes place in the algebra $\Q[p,q,r,z]$ of $\overline{H}$-invariants viewed
as a graded module over the algebra $\Q[a,b,c,d]$ of $\overline{G}$-invariants.
As a graded basis we use $p^i q^j r^k z^l$ with  $0 \leq i,j,k< 2$ and $0 \leq l < 30$.  Repeatedly using the vector equation in \S\ref{invariants2},
we expand the products
\[
\alpha_{e} p^i q^j r^k z^l= \sum_{I,J,K,L} M(e)^{i,j,k,l}_{I,J,K,L} p^I q^J r^K z^L
\]
to represent the covariants $\alpha_e$ as $240$-by-$240$ matrices $M(e)$
with entries in $\Q[a,b,c,d]$.   The general covariant 
\begin{equation}
\label{Zcovariant}
Z = s \alpha_1 +  t \alpha_7 + u \alpha_{13} + v \alpha_{19}
\end{equation}
satisfies the characteristic polynomial of $M = s M(1) + t M(7) + u M(13) + v M(19)$.
In other words, $Z$ satisfies a degree $240$ polynomial equation 
\[
F(A,B,C,D,Z) = Z^{240} + c_2 Z^{228} + c_3 Z^{222} + c_4 Z^{216} + c_5 Z^{210} + \cdots = 0
\]
with $F$ from \eqref{deg240}.  We need to calculate $A$, $B$, $C$, $D$ in terms of the free parameters
 $a$, $b$, $c$, $d$, $s$, $t$, $u$, and $v$.  

Define normalized traces $\tau_n$ by 
\[
6 \tau_n = \Tr(M^{6n}) = \sum_{i+j+k+l=6n} \binom{6n}{i,j,k,l} s^i t^j u^k v^l \Tr(M(1)^i M(7)^j M(13)^k M(19)^l).
\]
Because the first trace $\tau_1$ is $0$, standard symmetric polynomial formulas simplify, giving 
$(c_{2},c_3,c_4,c_5)=(-\tau_2/2,\tau_3/3,\tau_2^2/8-\tau_4/4,\tau_2 \tau_3/6 - \tau_5/5)$.  Then \eqref{deg240} yields 
 \begin{equation}
 \label{ABCDviatraces}
(A,B,C,D) = \left(
\frac{ {-\tau_2}}{30240},\frac{ {-\tau_3}}{7862400},\frac{3667
     {\tau^2_2}-5600  {\tau_4}}{9390915072000},\frac{  2521
     {\tau_2}  {\tau_3}-  2688  {\tau_5}}{886312627200000} 
     \right).
 \end{equation}
 
 The matrices $M^k$ have entries in $\Q[a,b,c,d,s,t,u,v]$ and for $k=1$, \dots, $6$ they take approximately 
 $2$, $10$, $40$, $125$, $300$, and $675$ megabytes to store.   The matrix $M^6$ suffices to determine
 $A$ because the evaluation of $\Tr(M^{12}) = \Tr(M^6 \cdot M^6)$ does not require the full matrix
 multiplication on the right.  However we would not be able to continue in this way to the needed $M^{15}$.  
In contrast, the $M(e)$ have entries only in $\Q[a,b,c,d]$ and take less space to store. 
 The worst of the $M(e)^j$ that we actually use in the 
 above expansion is  $M(19)^{15}$, which requires about $210$ megabytes to store.  
 By getting the terms in smaller batches and discarding matrix products when no longer needed, 
 we can completely compute all of $A$, $B$, $C$, and $D$ without memory overflow.
   In principle, one could repeat
 everything in the contravariant case, although here the initial matrix $M^*$ 
 takes twice as much space to store as $M$.

The polynomials $A$, $B$, $C$, and $D$
have respectively degrees $12$, $18$, $24$, and $30$ in $s$, $t$, $u$, and $v$.  
Also,  assign weights 
$(12, 18,24, 30, -1, -7, -13, -19)$ to $(a,b,c,d,s,t,u,v)$. 
Then all 
four polynomials are homogeneous of weight zero. 
The bigradation allows $A$, $B$, $C$, and $D$ to have
$14671$, $112933$, $515454$, and
$1727921$ terms respectively.  With our choice
of $\alpha_{13}$ and $\alpha_{19}$, respectively $67$, $170$, $100$, and $824$
of these terms vanish, so $A$, $B$, $C$, and $D$ have the
number of terms reported in the introduction.

Not only do the polynomials have many terms, but the coefficients can 
have moderately large numerators.   The largest absolute value of 
all the numerators is achieved by the term
\[
2^{30} \cdot 3^3 \cdot 5^{23}  \cdot 1381131815224116413 \cdot a^3  b  c^5  d^{10} u^{16} v^{14}
\]
in $D$.  On the another hand, denominators of the coefficients in $A$, $B$, $C$, and $D$
always divide $5$, $5^2$, $5^3$, and $5^5$ respectively.

\subsection{Geometric summary}
We now summarize our results in the following theorem. The~$\rhobar$ of \S1.2 is the mod 3 representation of the initial genus two curve \eqref{quintic1}.  So, to be more explicit, we write $\MM_{a,b,c,d} = \MM_2^w(\rhobar)$ below.

\begin{thm2}
\label{thm2}

Fix an equation $y^2=x^5 + a x^3 +b x^2 + c x + d$  
defining a curve $X$ over $\Q$.
Let $\MM_{a,b,c,d}$ be the moduli space of pairs $(Y,i)$ with $Y$ a Weierstrass curve  and 
$i : \Jac(X)[3] \rightarrow \Jac(Y)[3]$ a symplectic isomorphism on the $3$-torsion points of their
Jacobians.   Then $\MM_{a,b,c,d}$ can be realized as the complement of a discriminant locus
$\MZ_{a,b,c,d}$ in the projective three-space $\Proj \; \Q[s,t,u,v]$.   The covering maps to the moduli space $\MMW  \subset \Proj \; \Q[A,B,C,D]$
have degree $25920$ and are given by 
\begin{equation}
\label{maineq2}
(A,B,C,D)=(A(a,\dots,v),B(a,\dots,v),C(a,\dots,v),D(a,\dots,v)).
\end{equation}
The formula  
\begin{equation}
\label{ucurve2}
y^2 =  x^5 + A(a,\dots,v) x^3 + B(a,\dots,v) x^2 +  C(a,\dots,v) x +  D(a,\dots,v)
\end{equation}
gives the universal Weierstrass curve $X(s,t,u,v)$ over $\MM_{a,b,c,d}$. 
\end{thm2}

\noindent The discriminant locus
$\MZ_{a,b,c,d}$ is given by the vanishing of the discriminant 
\begin{equation}
\label{disc2}
\Delta(A(a,\dots,v),\dots,D(a,\dots,v)) = \Delta(a,b,c,d) \delta(a,b,c,d,s,t,u,v)^3.
\end{equation}
Here $\Delta(a,b,c,d)$ is a nonzero constant while
$\delta(a,b,c,d,s,t,u,v)$ is homogeneous of degree forty in $s$, $t$, $u$, $v$. Geometrically, $\MZ_{a,b,c,d}$ is the union of forty planes and these
planes are permuted by $\Gal(\overline{\Q}/\Q)$ according to the roots 
of $p_{40}$ from the end of \S\ref{3tors2}.  
In principle, Theorem~\ref{thm2}
has a direct analog for $\MM_{a,b,c,d}^* \rightarrow \MMW$.  The
computer file only gives the starred coefficients evaluated at 
$(a,b,c,d,1,0,0,0)$, as this is sufficient for moving from 
one moduli space to the other.

\subsection{Finding $(s,t,u,v)$}
\label{finding2}
       Let $X$ and $Y$ be Weierstrass curves over 
 $\Q$ having isomorphic $3$-torsion and given by coefficient sequences $(a,b,c,d)$ 
 and $(A,B,C,D)$ respectively.  Then finding associated rational $(s,t,u,v)$
 is both theoretically and computationally more complicated 
 than in the genus one case of \S\ref{finding1}.
 
       As in the genus one case, for \eqref{maineq2} to have a solution, the ratio
  $\Delta_X/\Delta_Y$ must be a perfect cube by \eqref{disc2}.   Similarly, for the starred
  version of \eqref{maineq2} to have a solution
  the product $\Delta_X \Delta_Y$ must be a  
  perfect cube.
      The theoretical complication was introduced at the end of \S\ref{3tors2}:  
      the class modulo  
  cubes of the discriminant now depends on the model via 
  $\Delta(u^4 A,u^6 B,u^8 C,u^{10} D) = u^{40} \Delta(A,B,C,D)$.
  So as a preparatory step one needs to adjust the model of 
  $Y$ to some new $(A,B,C,D)$ before seeking solutions to \eqref{maineq2},
  and also to some typically different $(A^*,B^*,C^*,D^*)$
  before seeking solutions to the starred analog of \eqref{maineq2}.

        Having presented $Y$ properly, one then encounters the computational problem.  Namely both 
 \eqref{maineq2} and its starred version have
  $155520$ solutions $(s,t,u,v) \in \C^4$, 
  and so one cannot expect to find the rational ones by algebraic 
  manipulations.   Working numerically instead, 
 one gets $240$ solutions 
 $(p,q,r,z) \in \C^4$ to the large vector equation in 
 \S\ref{invariants2}.  Eight of these
 solutions are in $\R^4$.   These vectors
 yield eight vectors $(\alpha_1,\alpha_7,\alpha_{13},\alpha_{19}) \in \R^4$ from 
 the covariants in \S\ref{covariants2}, 
 and also eight vectors $(\beta_{11},\beta_{17},\beta_{23},\beta_{29}) \in \R^4$.
 Let $Z$ and $Z^*$ respectively run over the eight real roots
 of $F(A,B,C,D,U)$ and $F(A^*,B^*,C^*,D^*,U)$.  
 Then one can apply the LLL algorithm to find low height relations
 of the form \eqref{Zcovariant} and its starred variant 
 \[
 Z^* = s \beta_{11}+t \beta_{17} + u \beta_{23} + v \beta_{29}.
 \]
 When the image of $\Gal(\overline{\Q}/\Q)$ on $3$-torsion is sufficiently
 large then there will just be a single pair of solutions $\pm (s,t,u,v)$ from
 the eight equations of one type and none from the other eight equations.  
  The computer file provides a {\em Mathematica} program \verb@findisos@
   to do all steps at once.   Examples are given in \S\ref{richelot} and \S\ref{modularity}.

\section{Complements}
\label{complements}
    The four subsections of this section can be read independently.  
    
 \subsection{A matricial identity}
 \label{matricial}
The polynomials $A$, $B$, $C$, and $D$ in Theorem~\ref{thm2} satisfy the matricial identity
\begingroup
{\footnotesize
\begin{equation*}
\label{eq:linear}
\begin{aligned}
& \mathcal{E}(A(a,\dots,v),B(a,\dots,v),C(a,\dots,v), D(a,\dots,v),S,T,U,V) =  \mathcal{E}(a,b,c,d,M (S, T, U, V)^t), \end{aligned}
\end{equation*}
}%
\endgroup
where $\mathcal{E}$ can be any one of $A$, $B$, $C$, $D$, and 
$M$ is a $4\times4$ matrix with entries in $\Q[a,b,c,d,s,t,u,v]$ whose first column is $(s, t, u, v)^t$. The columns of~$M$ 
are homogeneous of degrees $1,7,13,19$ in $s,t,u,v$, and the rows are homogeneous of degrees $-1,-7,-13,-19$ with respect to the weights assigned in \S\ref{newcoeff2}.  

The situation in the $g=1$ case is analogous but enormously simpler:
{\small
{\renewcommand{\arraycolsep}{2pt}
\[
\begin{array}{rclc}
\Aae(\Aae(\aae,\abe,s,t),\Abe(\aae,\abe,s,t),S,T) & = & \Aae(\aae,\abe,M (S,T)^t), & \multirow{2}{*}{$\;\;M \! = \!
\left( \begin{matrix} s & - \aae s^2 t - 3 \abe s t^2 +  \aae^2 t^3/3 \\ t & s^3 + \aae s t^2 + \abe t^3 \end{matrix} \right).$} \\
\Abe(\Aae(\aae,\abe,s,t),\Abe(\aae,\abe,s,t),S,T) & =  & \Abe(\aae,\abe,M (S,T)^t),  
\end{array}
\]
}}%
Here, as is visible, columns have degrees $1$ and $3$ in $s$, $t$, while rows have
weights $-1$ and $-3$ with respect to the weights assigned in \S\ref{newcoeff1}.   

The identities say that changing the initial Weierstrass curve to a different one in $\MM_{a,b,c,d}$ has the effect of changing the parametrization of the family through a linear   
transformation $M$ of the covariants. In fact, our first method of calculating the quantities $\mathcal{E}(a,\dots,v)$ exploited this ansatz. Starting from a few curves with $a = b = 0$,  computing covariants numerically, and changing bases so as to meet the bigradation conditions of \S\ref{newcoeffs2}, we obtained the polynomials $\mathcal{E}(0,0,c,d,s,t,u,v)$.  
We then examined the matricial identity with $a=b=0$.  Comparing certain monomial coefficients, we determined the second column of $M$ precisely, the third column up to one free parameter, and the fourth column up to two free parameters. This exactly corresponds to the ambiguity in the covariants in degrees $13$ and $19$ described in \S\ref{covariants2}. Once a choice of $M$ was made, comparing coefficients again and solving the resulting linear equations determined the polynomials $\mathcal{E}(a,\dots,v)$ completely.

\subsection{Examples involving Richelot isogenies}
\label{richelot}
      Let $X$ and $Y$ be Weierstrass curves and let $I: \Jac(X) \rightarrow \Jac(Y)$ be
   an isogeny with isotropic kernel of type $(m,m)$ with $m$ prime to $3$.  Then $I$ induces an isomorphism 
   $\iota : \Jac(X)[3] \rightarrow \Jac(Y)[3]$ which is symplectic if $m \equiv 1 \; (3)$
   and antisymplectic if $m \equiv 2 \; (3)$.   In the following examples, $m=2$.

   Let $X_{e,f,g}$ be
   defined by \eqref{quintic1} with $(a,b,c,d)=$
   \begin{equation*}
   \left( -5 (7 e^2-2 f),-10 e (3 e^2-2 f),5 (32
    e^4-39 e^2 f+g),-4 e (24 e^4+115 e^2 f-5 g) \right).
   \end{equation*}
   The discriminant of $X_{e,f,g}$ is 
   \begin{equation*}
   \Delta_X = -2^{12} 5^5 \left(125 e^4+20 f^2-4 g\right)^2 \left(25 e^2 f-g\right)
    \left(25 e^2 f+g\right)^2.
   \end{equation*}
   Define $Y_{e,f,g}$ to be the quadratic twist by $2$ of $X_{e,-f,g}$.  
   The form of $(a,b,c,d)$ has been chosen so 
   that there is a Richelot isogeny from $\Jac(X_{e,f,g})$ to 
   $\Jac(Y_{e,f,g})$.  
   
   Let $\bar{\cdot}$ be 
   the involution
    of $\Q[e,f,g]$ given by $(\bar{e},\bar{f},\bar{g}) = (e,-f,g)$.   
    To make $\Delta_X \Delta_Y$ a cube and avoid denominators in $(s,t,u,v)$, present $Y_{e,f,g}$ via $(A,B,C,D) = (\bar{a} z^2,\bar{b} z^3,\bar{c} z^4,\bar{d} z^5)$ with $z = 2^3 5^4 \left(125 e^4+20 f^2-4 g\right)^4 \left(25 e^2 f+g\right)^6$.
   Applying the numeric method of \S\ref{finding2} and interpolating strongly suggests
   $(s,t,u,v)=$
   \[
   \pm  \left(-4 e (80 e^4+7 e^2 f-g), 2(40 e^4-9 e^2 f-g),-4
    e (5 e^2+2 f),5 e^2+ 2 f\right).
   \]
   Specializing the contravariant matrix $M(a,b,c,d,s,t,u,v)^*$ of \S\ref{newcoeff2} to
   $M(e,f,g)^*$ allows direct computation
   of its powers up through the needed fifteenth power.   Applying \eqref{ABCDviatraces} indeed recovers $(A,B,C,D)$ so that the interpolation was correct. 
     
      The examples of this subsection are already much simpler than
   the general case with its millions of terms.   For a smaller family of
   even simpler examples, now with all mod $3$ representations non-surjective, one can set $e=0$.  Then $b$, $d$, $B$, $D$, $s$, and $u$ are all $0$,
   while $a$, $c$, $A$, $C$, $t$, and $v$ are given by tiny formulas.

 \subsection{Explicit families of modular abelian surfaces} 
 \label{modularity}
 Our main theorem gives a process by which modularity of
 a genus two curve can be transferred to modularity of 
 infinitely many other genus two curves:
 \begin{cor2}
 \label{cor1}
  Suppose the genus two curve $X : y^2 = x^5 + a x^3 + b x^2 + c x + d$ has
  good reduction at $3$, and assume that~$A = \Jac(X)$
  satisfies all the conditions of~\cite[Prop.~10.1.1, 10.1.3]{BCGP}, so that~$X$ is modular.
  Then all the curves $X(s,t,u,v)$ or $X^*(s,t,u,v)$ having good reduction at $3$ are also 
  modular. 
 \end{cor2}
 \noindent The conclusion follows simply because the hypotheses imply that
 the new Jacobians 
 also satisfy the conditions of~\cite[Prop.~10.1.1, 10.1.3]{BCGP} and are thus 
 modular.  In particular, for any $(s,t,u,v) \in \mathbf{P}^3(\Q)$ reducing 
  to~$(1,0,0,0) \in \mathbf{P}^3(\mathbf{F}_3)$, the curves $X$ and $X(s,t,u,v)$
  are identical modulo $3$ and therefore $X(s,t,u,v)$ is modular.   
  
  The hypotheses of \cite[Prop.~10.1.1, 10.1.3]{BCGP} include that 
  the mod $3$ representation $\rhobar$ is not surjective.  
  The easiest way to satisfy the hypotheses is to look among $X$
  for which the geometric endomorphism ring of $\mbox{Jac}(X)$ 
  is larger than $\Z$.  One such $X$, appearing in \cite[Example 3.3]{CGS}, is given by
  \[
  (a,b,c,d) = \left( \frac{12}{5}, \frac{12}{5^2}, \frac{292}{5^3}, - \frac{3672}{5^5} \right),
  \]
  having arisen from the simple equation 
  $y^2 = (x^2 + 2 x+2)(x^2+2)x$.   This curve has conductor 
  $2^{15}$ and discriminant $\Delta_X=2^{23}$.  
  Applying the corollary, one gets infinitely many modular
  genus two curves $X(s,t,u,v)$. 
  For generic parameters, the geometric endomorphism ring of 
 $\Jac(X(s,t,u,v))$ is just $\Z$.

   It is much harder to directly find curves $Y$ satisfying the hypotheses of
  \cite[Prop~10.1.1, 10.1.3]{BCGP} and also satisfying 
  $\mbox{End}(\Jac(Y)_{\overline{\Q}}) = \Z$.  A short list 
  was found in~\cite{CGS}. The curve $Y$ in Example~3.3 there
  has 
  \[
  (A,B,C,D) =  \left( 
  \frac{2^7}{5},\frac{2^{11} \cdot 57}{5^2},-\frac{2^{12} \cdot 503}{5^3},\frac{2^{17} \cdot 17943}{5^5} \right)
  \]
  and comes from the simple equation $y^2 = (2 x^4 + 2 x^2 + 1)(2x+3)$.   
  It has conductor $2^{15} 5$ and  Example~3.3 also observes that its $3$-torsion is
  isomorphic to that of $X$.  
  
  While $Y$ was found in \cite{CGS} via an {\em ad hoc} search, it now
  appears as just one point in an infinite family.  To see this explicitly, 
  note that $\Delta_Y =  2^{83} 5^6$ so that $\Delta_Y/\Delta_X$ 
  is a perfect cube.   Numerical computation as in \S\ref{finding2}
   followed by algebraic 
    verification yields
  \[
  Y = X\left(\frac{129}{125},\frac{11}{25},\frac{3}{100},\frac{1}{20}\right).
  \]
  If this procedure had failed, we would have found the proper $X^*(s,t,u,v)$ by 
  dividing $(A,B,C,D)$ by $(2^4,2^6,2^8,2^{10})$ to make $\Delta_X \Delta_Y$ a cube.

\subsection{Analogs for $p=2$}  
\label{analogs}
Complex reflection groups also let one respond 
to the problem of the introduction for residual prime $p=2$ 
and dimensions $g=2$, $3$, and $4$ via descriptions of moduli spaces
related to $\MA_g(\rhobar)$.    A conceptual simplification is that
since $p=2$ one does not have the second collection
of spaces $\MA^*_g(\rhobar)$.  Correspondingly, 
 the relevant groups are actually reflection groups defined
over $\Q$, so that covariants and contravariants coincide.  
The cases of dimension $g=3$, $4$ make fundamental
use of work of Shioda \cite{Shioda}.   

We begin with the easiest case $g=2$, because it shows
clearly that our approach has classical roots in
Tschirnhausen transformations.  
Greater generality would be possible by using the 
symmetric group $S_6$, but we describe things instead using 
$S_5$ to stay in the uniform context of Weierstrass curves.
Let $\alpha_1$ be a companion
matrix of $x^5 + a x^3 + b x^2 + c x + d$.  
For $j=2$, $3$, $4$, let $\alpha_j=\alpha_1^j-k I$
where $k$ is chosen to make $\alpha_j$ 
traceless.   Then the curve
\[
y^2 = \det(x I - s \alpha_1 - t \alpha_2 - u \alpha_3 - v \alpha_4)
\]
has the same $2$-torsion as the original curve.
From this fact follows a very direct
analog of Theorem~\ref{thm2}, with 
the new $\MM_{a,b,c,d} \subset \Proj \; \Q[s,t,u,v]$
now mapping to the same 
$\MMW \subset \Proj \; \Q[A,B,C,D]$
with degree $120$.  Carrying out this easy computation,
the elements $A$, $B$, $C$, and $D$ of 
$\Q[a,b,c,d,s,t,u,v]$ respectively have $24$, $86$, $235$,
and $535$ terms.   Of course there is
nothing special about degree $5$, 
and the analogous computations in 
degrees $2g+1$ and $2g+2$ give statements 
about genus $g$ hyperelliptic curves with fixed $2$-torsion.

For $g=3$, we work with the moduli space $\MM^q_{3}$ of smooth plane
quartics which maps isomorphically to an open subvariety of $\MA_3$.  
From the analog addressed in \cite{CC}, 
we suspect that the varieties $\MA_3(\rhobar)$
are in general not rational.  To place ourself in a 
clearly rational setting, we work with the moduli space 
$\MM_{3}^f$ of smooth plane quartics with a rational flex. 
This change is analogous to imposing a rational Weierstrass
point on a genus two curve, although now the resulting
cover $\MM^f_{3} \rightarrow \MM^q_{3}$ has degree
twenty-four.  A quartic curve with a rational flex can always
be given in affine coordinates by
 \begin{equation}
\label{e7}  y^3 + (x^3 + a_8 x + a_{12}) y +(a_2 x^4 + a_6 x^3 + a_{10} x^2 + a_{14} x + a_{18})  =   0.
 \end{equation}
 Here the flex in homogeneous coordinates is at $(x,y,z)=(0,1,0)$ and its tangent line is the line at 
 infinity $z=0$.   Changing $a_d$ to $u^d a_d$ gives an isomorphic curve via 
 $(x,y) \mapsto (u^4 x,u^6 y)$.   The variety $\MM^f_{3}$ is
 the complement of a discriminant locus in the weighted
 projective space 
 $\Proj \; \Q[a_2,\dots,a_{18}] = \bbP^6(2,\dots,18)$.
 The invariant theory of the reflection group 
 $\ST36 = W(E_7) = C_2 \times \Sp_6(\F_2)$ 
 gives polynomials $A_i(a_2,\dots,a_{18},s_{-1},\dots,s_{-17})$ 
 of degree $i$ in the $s_{-j}$ and total weight $0$.  Following the template of
 the previous cases, for fixed $(a_2,\dots,a_{18})$ one has a six-dimensional variety 
 $\MM_{a_2,\dots,a_{18}} \subset \Proj \; \Q[s_{-1},\dots,s_{-17}]$
 parametrizing genus three curves with a rational flex and $2$-torsion 
 identified with that of \eqref{e7}.
 The covering maps $\MM_{a_2,\dots,a_{18}} \rightarrow \MM^f_{3}$ 
 now have degree $|\Sp_6(\F_2)| = 1451520$.   The number of terms 
 allowed in $A_i(a_2,\dots,a_{18},s_{-1},\dots,s_{-17})$ by the bigradation
 is the coefficient of $x^i t^{19 i}$ in 
 \begin{equation}
 \label{genfunct}
 \prod_{d \in \{2,6,8,10,12,14,18\}} \frac{1}{(1-t^d)(1-x t^d)}.
 \end{equation}
 For $i = 18$, this number is $11,617,543,745$, so complete computations
 in the style of this paper seem infeasible.  
 
 For $g=4$, one needs to go quite far away from the $10$-dimensional
 variety $\MA_4$ to obtain a statement parallel to the
 previous ones.  Even the nine-dimensional variety $\MM_4$ is too
 large because for a generic genus four curve $X$ corresponding to 
 a point in $\MM_4$, the image of $\Gal(\Qbar/\Q)$ in its action
 on $\Jac(X)[2]$ is $\Sp_8(\F_2)$, and this group is not
 a complex reflection group.  However, one can work with the smooth curves
  \begin{equation}
  \label{e8}
  y^3 + (a_2 x^3 + a_8 x^2 + a_{14} x + a_{20}) y + (x^5 + a_{12} x^3 + a_{18} x^2 + a_{24} x + a_{30})  =  0 
 \end{equation}
 and a corresponding seven-dimensional moduli space 
 $\MM_4^s \subset \bbP^7(2,\dots,30)$.   For a generic curve in \eqref{e8}, the
 image of $\Gal(\Qbar/\Q)$ is the index 136 subgroup 
 $\OO_8^+(\F_2):2$ of $\Sp_8(\F_2)$.  Now from the invariant theory of 
 the largest Shephard--Todd group $\ST37 = W(E_8) = 2.\OO_8^+(\F_2) \td 2$,
 one gets polynomials $A_i(a_2,\dots,a_{30},s_{-1},\dots,s_{-29})$
 and covering maps $\MM_{a_2,\dots,a_{30}} \rightarrow \MM^s_{4}$ 
 of degree $|\OO_8^+(\F_2) \td 2|=348,364,800$.   Aspects of 
 this situation are within computational reach; for example 
 Shioda computed the degree $240$ polynomial $F(a_2,\dots,a_{30},z)$ 
 analogous to \eqref{division3} and \eqref{deg240}.  However the number of allowed terms 
 in $A_i(a_2,\dots,a_{30},s_{-1},\dots,s_{-29})$ is even 
 larger than in the previous $g=3$ case, being  the coefficient of $x^i t^{31 i}$ in the
 analog of \eqref{genfunct} where $d$ runs over $\{2, 8, 12, 14,18,20,24,30\}$.  
 For $i=30$, this number is $100, 315, 853, 630, 512$.   
 We close the paper with this  $W(E_8)$ case because
 it is here that the paper actually began:  the polynomial \eqref{deg240} for 
 our main case $C_3 \times \Sp_4(\F_3)$ is
 also the specialization $F(0,0,a_{12},0,a_{18},0,a_{24},a_{30},z)$ 
 of Shioda's polynomial.

 \bibliographystyle{alpha}
\bibliography{Fixed3TorsionA}

\end{document}